\documentclass[10pt,a4paper,reqno]{amsart}

\usepackage[english]{babel}
\usepackage[latin1]{inputenc}
\usepackage[OT1]{fontenc} 

\usepackage{amssymb}

\def\qbino#1,#2,{{#1\atopwithdelims[]#2}}

\newcommand{\pref}[1]{{(\protect\ref{#1})}}

\newcommand\mc{\mathcal}
\def\newop#1{\expandafter\def\csname #1\endcsname{\mathop{\rm #1}\nolimits}}
\def\disp{\displaystyle}

\def\emm#1,{{\em#1}}

\newcommand\blsk{\baselineskip}

\newcommand\cd{\cdot}

\newcommand\du{\mathop{{\cup}\kern -.48em {\cdot}\kern 3 pt}}

\def\ch#1,#2,{{#1\choose #2}}

\theoremstyle{plain}
\newtheorem{prop}{Proposition} 
\newtheorem{lemma}[prop]{Lemma}
 
\newtheorem{theorem}[prop]{Theorem}

\newtheorem{proposition}[prop]{Proposition}
\newtheorem{corollary}[prop]{Corollary}
\newtheorem{conjecture}[prop]{Conjecture}

\newtheorem{definition}[prop]{Definition}

\newcommand\st{\; | \;}

\newcommand\cl[1]{{\mc #1}}

\newcommand\cls{{\mc S}}

\newcommand\des{\mathop{\mathrm{des}}}
\def\newMAH#1{%
\expandafter\def\csname #1\endcsname{\mathop{\mbox{{\sc#1}}}\nolimits}%
}
\newMAH{maj}
\newMAH{mil}
\newMAH{inv}
\newMAH{mak}
\newMAH{stat}
\newcommand\lmak{\ell\mak}

\def\newexpMAH#1{%
\expandafter\def\csname exp#1\endcsname{\mathop{\mbox{{\footnotesize\sc{#1}}}}\nolimits}%
}

\newexpMAH{inv}
\newexpMAH{maj}
\newexpMAH{mil}

\newop{ros}
\newop{rcs}
\newop{rcb}
\newop{rob}
\newop{los}
\newop{lcs}
\newop{lcb}
\newop{lob}
\newop{lsb}
\newop{rsb}

\newop{scl}

\newMAH{ROS}
\newMAH{RCS}
\newMAH{RCB}
\newMAH{ROB}
\newMAH{LOS}
\newMAH{LCS}
\newMAH{LCB}
\newMAH{LOB}
\newMAH{LSB}
\newMAH{RSB}

\newexpMAH{ros}
\newexpMAH{rcs}
\newexpMAH{rcb}
\newexpMAH{rob}
\newexpMAH{los}
\newexpMAH{lcs}
\newexpMAH{lcb}
\newexpMAH{lob}
\newexpMAH{lsb}
\newexpMAH{rsb}

\newop{open}
\newop{clos}
\newop{block}
\newcommand\blk{\block}

\def\bmaj{\mathop{\mbox{\hbox{\rm b}\hbox{$\maj$}}}}
\def\cbmaj{\mathop{\mbox{\hbox{\rm cb}\hbox{$\maj$}}}}

\def\binv{\mathop{\mbox{\hbox{\rm b}\hbox{$\inv$}}}}
\def\cbinv{\mathop{\mbox{\hbox{\rm cb}\hbox{$\inv$}}}}

\def\cmajlsb{\mathop{\mbox{\hbox{\rm cmaj}\hbox{$\LSB$}}}}
\def\cinvlsb{\mathop{\mbox{\hbox{\rm cinv}\hbox{$\LSB$}}}}

\def\bmajmil{\mathop{\mbox{\hbox{\rm bmaj}\hbox{$\mil$}}}}

\newop{ncl}

\def\sq#1,#2,{\mbox{$S_q(#1,#2)$}}
\def\sqt#1,#2,{\mbox{$\tilde{S}_q(#1,#2)$}}
\def\aq#1,#2,{\mbox{$A_q(#1,#2)$}}

\def\qbino#1,#2,{{#1\atopwithdelims[]#2}}

\def\p#1,#2,{\ensuremath{\cl{P}_{#1}^{#2}}}

\def\op#1,#2,{\ensuremath{\cl{O}\cl{P}_{#1}^{#2}}}

\begin{document}

\thispagestyle{empty}

\begin{center} {\Large \bf Statistics on ordered partitions of
    sets}\\

\renewcommand{\thefootnote}{\fnsymbol{footnote}}

\vspace{12pt} {\large Einar Steingr\'{\i}msson\footnote{Part of this
    research was carried out while the author was visiting Université
    Louis-Pasteur in Strasbourg in 1996, supported by the EU Network
    in Algebraic Combinatorics.}\\
  {Department of Computer and Information Sciences\\ University of Strathclyde\\Glasgow, UK}\\
{\normalsize\tt einar@alum.mit.edu}}

\end{center}

\begin{abstract}
  We introduce several statistics on ordered partitions of sets, that is, set partitions where the blocks are permuted arbitrarily.  The distribution of these statistics is closely related to the $q$-Stirling numbers of the second kind.  Some of the statistics are generalizations of known statistics on set partitions, but others are entirely new.  All the new ones are sums of two statistics, inspired by statistics on permutations, where one of the two statistics is based on a certain partial ordering of the blocks of a partition.
\end{abstract}




\begin{abstract}
  We introduce several statistics on ordered partitions of sets, that
  is, set partitions where the blocks are permuted arbitrarily.  The
  distribution of these statistics is closely related to the
  $q$-Stirling numbers of the second kind.  Some of the statistics are
  generalizations of known statistics on set partitions, but others
  are entirely new.  All the new ones are sums of two statistics,
  inspired by statistics on permutations, where one of the two
  statistics is based on a certain partial ordering of the blocks of a
  partition.
\end{abstract}

\emm Keywords:, Ordered set partitions, q-Stirling numbers, permutation
statistics.

\setcounter{section}{-1} 

\section{Prologue}

This paper was first made public on my webpage at Chalmers University of Technology in 2001, after it had been rejected by two journals. It was first posted on the arXiv in 2006 \cite{es-stat-ord-part}, with (trivial) updates in 2007 and 2014.  

The reason it is being published now is to put it more formally on the record, as it has been cited in several papers, dealing both directly with some of the conjectures in the paper and in work in other areas. In a series of papers, Ksavrelof and Zeng \cite{ksavrelof-zeng-nouvelle-stat}, Ishikawa, Kasraoui and Zeng \cite{ishikawa-kasraoui-zeng-em-stats-ordsetpart, ishikawa-kasraoui-zeng-survey} and Kasraoui and Zeng \cite{kasraoui-zeng-em-stats-ordsetpart-2} proved all the conjectures in this paper, the last of these papers giving combinatorial proofs.  Also, Remmel and Wilson \cite{remmel-wilson-macmahon-ordsetpart} solved the problem that was the original motivation for the present paper, namely finding a bijective proof of Proposition~\ref{zz} here, due to Zeng and Zhang~\cite{zeng-zhang-newton}, thus ``explaining'' combinatorially this strong connection between ordered set partitions and permutations. Wilson then extended this correspondence to multiset partitions and \emph{descent-starred permutations}~\cite{wilson-macmahon-ordmultisetpart}.  

The statistic $\ROS$ defined here also has some bearing on the Delta Conjecture of Haglund, Remmel and Wilson \cite{haglund-remmel-wilson-delta}, and played a role in Haglund, Rhoades and Shimozono's paper on that conjecture and generalized coinvariant algebras~\cite{haglund-rhoades-shimozono-ordsetpart-delta}.

Apart from this Prologue, the deletion of a now obsolete footnote and the updating of references, the content of the paper is the same as in the original version.

\section{Introduction}

The Stirling numbers of the second kind, $S(n,k)$, which count the
partitions of an $n$-element set into $k$ blocks, have been much
studied.  Their $q$-analog $\sq n,k,$, the $q$-Stirling numbers of the
second kind, can be defined by $\sq n,k, = 0$ if $k>n$ or $k<0$ and,
for $n\ge k\ge0$, by the identity
\begin{eqnarray}\label{sqdef}%
[k]!\cd\sq n,k, = \sum_i{(-1)^i\cd\qbino k,i,\cd q^{\ch i,2,}\cd [k-i]^n},
\end{eqnarray} 
where
$$
[k] = 1+q+\cdots+q^{k-1}, ~~~ [k]! = [k][k-1]\cdots[1], ~~\mbox{ and
  } ~~ \qbino n,k, = \frac{[n]!}{[k]!\cd[n-k]!}.
$$
Here, as in the remainder of the paper, sums are taken to be over all
integers, unless explicitly stated otherwise.

There is an alternative definition of the $q$-Stirling numbers.  Those
numbers, sometimes denoted $\sqt n,k,$, are related to $\sq n,k,$ by
$\sq n,k, = q^{\ch k,2,}\sqt n,k,$.  In fact, the factor $q^{\ch
  k,2,}$ is explicit in some of our statistics, indicating that they
essentially correspond to $\sqt n,k,$, but as this factor does not
appear transparently in some other statistics we take \pref{sqdef}
above as our definition.

There are several ways to define $\sq n,k,$ combinatorially, most of
them based on statistics on set partitions.  Perhaps the simplest of
these statistics (in terms of definition) is the one due to Milne
\cite{milne-q-rgf}.  It can be defined as follows: Given a partition of the
set $\{1,2,\ldots,n\}$ into $k$ blocks, order the blocks in the
standard way, that is, by increasing least element.  Let $b_i$ be the
size of the $i$th block in this ordering.  Milne's statistic, which we
call $\mil$, is then defined to be $b_2+2b_3+\cdots+(k-1)b_k$.  As an
example, $\mil(14-238-5-67)= 1\cd 3 + 2\cd 1 +3\cd 2 = 11$.

The following lemma is easily proved from identity \pref{sqdef} or
from the definition of $\mil$ given above.
\begin{lemma}\label{rec-lemma}%
  The numbers $\sq n,k,$ satisfy the recurrence
\begin{eqnarray}\label{sqrec}%
\sq n,k, = q^{k-1}\cd \sq n-1,k-1, + [k]\cd\sq n-1,k,.
\end{eqnarray} 
\end{lemma}

Many authors have studied statistics on set partitions with
distribution given by the numbers \sq n,k,.  Apart from the paper by
Milne \cite{milne-q-rgf}, who seems to have pioneered the study of partition
statistics whose distribution is given by the $q$-Stirling numbers, we
mention \cite{ER-juggling, garsia-remmel-q-rook, sagan-maj-stat, wachs-white-stirling,
  white-interpolating}.  Of these, the paper by Wachs and White
\cite{wachs-white-stirling} is the most comprehensive.  Also, there is a
substantial literature on various statistics on set partitions
restricted to so-called \emm non-crossing, partitions.  See for
example \cite{bona-simion, reiner-noncrossing, simion-noncross-partitions, white-interpolating}.

In the present paper we study statistics on \emm ordered, partitions
of sets, that is, set partitions where the blocks of a partition are
ordered arbitrarily.  The original motivation of this study is an
identity relating $q$-Stirling numbers of the second kind to
$q$-Eulerian numbers (see Proposition~\ref{zz} below).  The Eulerian
number $A(n,k)$ counts permutations in the symmetric group $\cls_n$
with $k$ \emm descents,.  A descent in a permutation $p = a_1a_2\cdots
a_n$ is an $i$ such that $a_i>a_{i+1}$.  There is a basic identity
relating the $S(n,k)$'s and the $A(n,k)$'s, namely
\begin{eqnarray}\label{stirl-to-eul}
k!\cd S(n,k)= \sum_i {\ch n-i,k-i,\cd A(n,i-1),}
\end{eqnarray}
which is easily proved combinatorially.  The $q$-analog of the
Eulerian numbers is the $\maj$-statistic, defined as the sum of the
descents in a permutation. As an example, the permutation
$p=5261743$ has descents 1, 3, 5 and 6, and thus
$\maj p=1+3+5+6=15$.  We denote by  $\aq n,k,$  the distribution of
permutations with $k$ descents according to $\maj$, that is, 
$$
\aq n,k, = \sum_{p\in\cls_n^k}{q^{\expmaj p}},
$$
where $\cls_n^k$ is the set of permutations in $\cls_n$ with exactly
$k$ descents.

The $q$-analog of identity \pref{stirl-to-eul}
was derived by Zeng and Zhang \cite{zeng-zhang-newton}, using analytic methods,
and goes as follows.
\begin{proposition}[Zeng and Zhang\protect{\cite[Proposition
    4.6]{zeng-zhang-newton}}]\label{zz}%
  For all $n$ and~$k$ with $0\le k\le n$ we have
\begin{eqnarray}\label{qsa}
[k]!\cd\sq n,k,= \sum_i {q^{k(k-i)}\cd\qbino n-i,k-i,\cd\aq n,i-1,}.
\end{eqnarray}
\end{proposition}
The distribution of the bistatistic $(\des,\maj)$ on the symmetric
group $\cls_n$, that is, the sum of $\aq n,k,$ over all $k$, has been
thoroughly studied.  This statistic, together with all other
bistatistics with the same distribution, is said to be \emm
Euler-Mahonian, (the number of descents is an \emm Eulerian, statistic
and $\maj$ is a \emm Mahonian, statistic).  Because of the relation
between the $q$-Eulerian numbers and the $q$-Stirling numbers encoded
by identity \pref{qsa}, we will refer to statistics on ordered
partitions of sets with the distribution $[k]!\cd\sq n,k,$ as
Euler-Mahonian, and we will even call Euler-Mahonian those statistics
on (unordered) partitions of sets that have the distribution $\sq
n,k,$.

As mentioned above, the original motivation of this paper was the
identity~\pref{qsa}.  We have found several statistics, defined on
ordered partitions of sets, whose distribution we conjecture to be
that given by either side of \pref{qsa}.  Two of these we prove to
have this distribution, but in neither case have we been able to find
a bijective proof.  The obvious combinatorial proof of identity
\pref{stirl-to-eul} does not generalize to the $q$-analog case of
identity \pref{qsa} in any straightforward way using the known
statistics on set partitions.  Thus, identity~\pref{qsa} still lacks a
combinatorial proof.

Now, given any statistic $S$ with the distribution $\sq n,k,$ on $\p
n,k,$ we could define a composite statistic $\stat=S+T$ on $\op n,k,$
by computing $S$ on the standard ordering of the blocks of $\pi\in \op
n,k,$ and letting $T$ be any Mahonian \emm permutation statistic,
(such as $\inv$, the number of inversions, or $\maj$) computed on the
permutation induced by the ordering of the blocks of $\pi$, since the
Mahonian permutation statistics have distribution $[k]!$ on $\cls_k$.
Thus, we would be computing two separate statistics and effectively
ignoring the ordering of the blocks for one of the statistics.
However, we know of no such statistic that lends itself to a
combinatorial proof of identity \pref{qsa}.

In fact, we are aware of only one statistic defined on set partitions
that is independent of the ordering of the blocks, and thus can be
combined with any Mahonian permutation statistic on the permutation of
the blocks to obtain a statistic with distribution $[k]!\sq n,k,$ on
$\op n,k,$.  This is the \emm intertwining number, of Ehrenborg and
Readdy \cite[\S 6]{ER-juggling}, an interesting statistic with
properties quite different from those considered here.  But, since
we don't make any use of this statistic we omit its definition.

Although we have not come up with a bijective proof of \pref{qsa} we
have found a statistic, which we call $\bmajmil$ (see
Definition~\ref{bmm-def}), defined on ordered partitions, whose
distribution is given by $[k]!\cd\sq n,k,$.  However, our proof
consists of constructing a bijection between ordered partitions of
$\{1,2,\cdots,n\}$ with $k$ blocks and the set of permutations with at
most $n-k$ descents, marked in certain ways, showing that the
statistic in question has the distribution given by the right hand
side of \pref{qsa}.  We also show that a straightforward
generalization of one of the statistics by Wachs and White
\cite{wachs-white-stirling} has distribution $[k]!\cd\sq n,k,$.

We have also found several other statistics that we conjecture to
have the distribution $[k]!\cd\sq n,k,$.  All of these conjectures
have been verified by computer for $n\le11$, so it seems unlikely that
they could be wrong.

It is somewhat intriguing that the new statistics introduced are all
sums of two statistics, where one is essentially one of two Mahonian
permutation statistics, but defined on a certain partial ordering of
the blocks of a partition.

It should also be mentioned that the ordered partitions studied here
have recently been treated by Krob, Latapy, Novelli, Phan and Schwer
\cite{krob-pseudo}, under the name of \emm pseudo-permutations,.  Their
results are quite different from ours, but it seems likely that some
connections will emerge.

\section{Generalizations of known statistics}

Let $\p n,k,$ be the set of (unordered) partitions of the set
$\{1,2,\ldots,n\}$ with $k$ blocks and let $\op n,k,$ be the set of
ordered partitions of the set $\{1,2,\ldots,n\}$ with $k$ blocks.  We
now introduce several statistics on $\op n,k,$ that are
generalizations of known statistics on $\p n,k,$.

So far, most known statistics on $\p n,k,$ can be defined in terms of
inversions between the letters (integers) in a partition $\pi$ and the
\emm openers, and \emm closers, of the blocks of $\pi$.  The opener of
a block is its least element and the closer is its greatest element.
For example, the partition $\pi=136-27-4-58$ has openers 1, 2, 4 and 5
and closers 6, 7, 4 and 8.

On unordered partitions we always assume, when referring to the
ordering of the blocks of a partition, that the blocks are written in
increasing order of their respective openers.  We call this the \emm
standard ordering,.

One statistic whose distribution on $\p n,k,$ is $\sqt n,k,$ is the
sum $\ros\pi=\sum_i{\ros_i\pi}$ where $\ros_i\pi$ is the number of
openers in $\pi$ that are smaller than $i$ and that belong to blocks
to the right of the block containing $i$ (the name $\ros$ is of course
an abbreviation of ``right opener smaller'').  For example, the values
of $\ros_i$ for the partition $\pi=136-27-4-58$ are 0, 1, 3, 0, 2, 0,
0, 0 (written in the order in which the letters appear in $\pi$), so
$\ros\pi=6$.

We now define ten partition statistics.  Four of these were defined by
Wachs and White \cite{wachs-white-stirling}, although their treatment was in
terms of restricted growth functions, a different way of representing
partitions.  Another four of our statistics are mirror images of the
aforementioned ones that contribute nothing new in the case of
unordered partitions.  The last two statistics, essentially defined by
Foata and Zeilberger \cite{foata-zeilberger-denerts} for permutations, are in fact each
equal to the difference of two of the first eight statistics.
However, they are useful to define and they also make clear the
similarity between some of our partition statistics and known
permutation statistics.

\enlargethispage{1\blsk}
\begin{definition}\label{stats-def}%
  Given a partition $\pi\in\op n,k,$, let $\open\pi$ and $\clos\pi$ be
  the set of openers and closers, respectively, of $\pi$.  Let
  $\blk(i)$ be the number (counting from the left) of the block
  containing the letter $i$.  We define eight \emm coordinate
  statistics, as follows:
\begin{eqnarray*}
&\ros_i \pi &=~~ \#\{j\st  i>j, ~~ j\in\open\pi, ~ \blk(j)>\blk(i)\},\\
&\rob_i \pi &=~~ \#\{j\st  i<j, ~~ j\in\open\pi, ~ \blk(j)>\blk(i)\},\\
&\rcs_i \pi &=~~ \#\{j\st   i>j, ~~ j\in\clos\pi, ~~\blk(j)>\blk(i)\},\\
&\rcb_i \pi &=~~ \#\{j\st   i<j, ~~ j\in\clos\pi, ~~\blk(j)>\blk(i)\},\\
&\los_i \pi &=~~ \#\{j\st  i>j, ~~ j\in\open\pi, ~ \blk(j)<\blk(i)\},\\
&\lob_i \pi &=~~ \#\{j\st  i<j, ~~ j\in\open\pi, ~ \blk(j)<\blk(i)\},\\
&\lcs_i \pi &=~~ \#\{j\st   i>j, ~~ j\in\clos\pi, ~~\blk(j)<\blk(i)\},\\
&\lcb_i \pi &=~~ \#\{j\st   i<j, ~~ j\in\clos\pi, ~~\blk(j)<\blk(i)\}.
\end{eqnarray*} 
Moreover, we let $\rsb_i$ be the number of blocks $B$ to the right of
the block containing $i$ such that the opener of $B$ is smaller than
$i$ and the closer of $B$ is greater than $i$ ($\rsb$ is an
abbreviation for ``right, smaller, bigger'').  Also, we define
$\lsb_i$ in an analogous way, with ``right'' replaced by ``left.''

We then set
$$
\ros\pi=\sum_i{\ros_i\pi}
$$
and likewise for the remaining nine statistics, i.e. each of
$\rob$, $\rcs$, $\rcb$, $\los$, $\lob$, $\lcs$, $\lcb$, $\rsb$, $\lsb$
is defined to be the sum over all $i$ of the respective coordinate
statistic.
\end{definition} 
Note that the statistic $\rsb$ was first defined, essentially, in
\cite{foata-zeilberger-denerts} (see also \cite{csz}) as a ``partial statistic'' in the
definition of the permutation statistic $\mak$.  We will treat a
partition statistic analogous to $\mak$ later in this paper.

As an example of the coordinate statistics just defined we give the
values of three of them on the ordered partition $\pi ~ =47-3-159-68-2$:
\begin{eqnarray*}%
\begin{array}{lcccccccccccccc}
\pi =  && 4 & 7 & - & 3 & - & 1 & 5 & 9 & - & 6 & 8 & - & 2\\[2pt]
\ros_i: && 3 & 4 &   & 2 &   & 0 & 1 & 2 &   & 1 & 1 &   & 0\cr
\lcb_i: && 0 & 0 &   & 1 &   & 2 & 1 & 0 &   & 2 & 1 &   & 4\\
\rsb_i: && 1 & 2 &   & 1 &   & 0 & 0 & 0 &   & 0 & 0 &   &  0\\
\end{array}
\end{eqnarray*} 

On unordered partitions, the statistic $\mil$ mentioned above is in
fact equal to $\los$ (left opener smaller), since the opener of every
block to the left of the block containing a given letter is smaller
than that letter.  Thus, $\los$ has the distribution $\sq n,k,$ on $\p
n,k,$.  On the other hand, $\lob$ (left opener bigger) is clearly
identically zero on $\p n,k,$.

On ordered partitions the situation is much simpler, since ``left''
and ``right'' have equal status, that is, reversing the order of the
blocks in an ordered partition turns a left opener into a right opener
and likewise for closers.  Thus every right statistic is
equidistributed with its corresponding left statistic.

Moreover, given a partition $\pi$, let $\pi^c$ be the partition
obtained by \emm complementing, each of the letters in $\pi$, that is,
by replacing the letter $i$ by $n+1-i$.  It is then easily checked
that $\rcb\pi^c =\ros\pi$ and that $\rcs\pi^c =\rob\pi$.  Thus the
eight statistics obtained by independently varying left/right,
opener/closer and smaller/bigger fall into only two categories when it
comes to their distribution on ordered partitions.  One of these
categories consists of $\ros$, $\rcb$, $\los$ and $\lcb$, and the
other contains $\rob$, $\rcs$, $\lcs$ and $\lob$.  Since
$\rcb_i\ge\rob_i$ for all $i$ and any partition $\pi$, it is clear
that the statistics in the latter category are ``smaller'' than those
in the first one.  Indeed, the ``small'' statistics do not have a
distribution related in any obvious way to the $\sq n,k,$, except that
$\rob_i$ equals $\rcb_i-\rsb_i$, which is easily proved as are the
analogous identities for $\rcs$, $\lcs$ and $\lob$.  However, the
``small'' statistics play a role in some of our new statistics, which
consist of combinations of these with yet other statistics.

It should perhaps be mentioned that there is a way to redefine $\ros$
in order to ``split'' it into a partition statistic and a Mahonian
statistic on the permutation of the blocks.  Namely, $\ros_i$, when
restricted to the openers of $\pi$, records just the inversions among
the openers, and thus the inversion statistic on the permutation of
the blocks of $\pi$ (as compared to the standard ordering).  However,
as mentioned before, we have not been able to exploit this to find
a combinatorial proof of Proposition~\ref{zz}.

The statistic $\ros$ has distribution $\sqt n,k,$ on $\p n,k,$ and on
$\op n,k,$ its distribution is $[k]!\sqt n,k,$.  In order to get the
distribution $[k]!\sq n,k,$ we therefore add $\ch k,2,$, where $k$ is
the number of blocks in $\pi$.  We call this new statistic $\ROS$, so
$\ROS\pi=\ros\pi+\ch k,2,$.

\begin{theorem}\label{ros-thm}%
  The statistic $\ROS$ is Euler-Mahonian on ordered partitions, that
  is,
$$
\sum_{\pi\in\op n,k,}{q^{\expros\pi}} = [k]!\;\sq n,k,.
$$
\end{theorem}
\begin{proof}{%
  The proof is by induction on $n$, with $k$ fixed.  Suppose first
  that $\pi$ is a partition in \op n,k,, such that $n$ is a singleton
  block.  Then $\pi$ can be obtained from a unique partition $\pi'$ in
  \op n-1,k-1,, by inserting $n$ as a singleton block in one of the
  slots between blocks in $\pi'$, or before, or after, all the blocks
  in $\pi$.  If $\pi$ is obtained by inserting $n$ right after block
  $i$ in $\pi'$, then there will be $(k-1-i)$ blocks following $n$ so
  this increases $\ROS$ by $k-1-i$.  Moreover, since we now have $k$
  blocks instead of $k-1$, $\ROS$ is additionally increased by $\ch
  k,2,-\ch k-1,2, = k-1$.  Therefore the total increase in $\ROS$ is
  $(k-1-i) + (k-1)$.

Thus, if $S$ is the set of partitions obtained from the partition
$\pi'$ by inserting $n$ in the $k$ slots between blocks in $\pi'$ or
before or after all of the blocks, then we have
$$
\sum_{\pi\in S}{q^{\expros\pi}} =
q^{\expros\pi'}q^{k-1}(1+q+\cdots+q^{k-1}) = q^{\expros\pi'}q^{k-1}[k].
$$

On the other hand, if $n$ is not a singleton block in $\pi$, then
$\pi$ can be obtained from a unique partition $\pi'$ in \op n-1,k,, by
inserting $n$ in one of the blocks in $\pi'$.  If $n$ is thus inserted
into the $i$-th block we increase $\ROS$ by $(k-i)$, so the
corresponding sum, with $T$ defined analogously to $S$ above, is
$$
\sum_{\pi\in T}{q^{\expros\pi}} =
q^{\expros\pi'}(1+q+\cdots+q^{k-1}) = q^{\expros\pi'}[k].
$$
Thus we have, assuming the statement true for $n-1$, that
\begin{eqnarray*}
\sum_{\pi\in\op n,k,}{q^{\expros\pi}} &=&
 q^{k-1}[k][k-1]!\sq n-1,k-1, ~+~ [k][k]!\sq n-1,k,\\
&=&  [k]!\;\biggl(  q^{k-1}\sq n-1,k-1, ~+~[k]\sq n-1,k,\biggr),
\end{eqnarray*}
which, by Lemma~\ref{rec-lemma}, equals $[k]!\;\sq n,k,$ as desired.
The case $n=1$ is trivial (and the case $n=0$ a matter of definition).
}\end{proof}

\begin{corollary}\label{etc-coro}%
The following statistics are Euler-Mahonian on $\op n,k,$:
\begin{eqnarray*}
\RCB&=&\rcb+\ch k,2,,\\
\LOS&=&\los\,+\,\ch k,2,,\\
\LCB&=&\lcb\,+\,\ch k,2,.
\end{eqnarray*} 
\end{corollary}

\section{A new Euler-Mahonian statistic on ordered  partitions}

We now define a new statistic that has the distribution $[k]!\sq n,k,$
on $\op n,k,$.  This statistic, and each of the new statistics we
introduce in the next section, is a combination of two statistics $S$
and $T$, where $S$ is a partition statistic that is inspired by a
permutation statistic.  What is surprising is that the $T$-statistics
are permutation statistics computed on the \emm partial,\/ ordering of
the blocks of a partition defined by setting $B_i < B_j$ if \emm
each,\/ letter of $B_i$ is smaller than \emm every,\/ letter of $B_j$.
We consider two such permutation statistics on the blocks,
corresponding to $\inv$ and $\maj$, respectively.
\begin{definition}\upshape%
  Let $\pi$ be an ordered partition in $\op n,k,$ with blocks
  $B_1,B_2,\ldots,B_k$ and let the partial ordering of the blocks be
  as
  defined in the preceding paragraph.\\[.5ex]
  We say that $i$ is a \emm block descent, \/ in $\pi$ if
  $B_i>B_{i+1}$.\\[.5ex]
  The \emm block major index, of $\pi$, denoted $\bmaj\pi$, is the sum
  of the block descents in $\pi$.\\[.5ex]
  A \emm block inversion, in $\pi$ is a pair $(i,j)$ such that $i<j$
  and $B_i>B_j$.\\[.5ex]
  The \emm block inversion number of $\pi$, denoted $\binv\pi$, is the
  number of block inversions in $\pi$.\\[.5ex]
  Moreover, we set $\cbmaj=\ch k,2,-\bmaj$ and $\cbinv=\ch
  k,2,-\binv$.
\end{definition} 
Note that $\cbmaj$ is the sum of the elements in the complement of the
set of block descents in $\pi$, whence the prefixed ``c,'' and
likewise for $\cbinv$.

For example, if $\pi= 41-96-5-87-32$ (where we write the elements of
each block decreasingly to emphasize the block descents) then $\binv\pi
= 0+2+1+1+0=4$ and $\bmaj\pi=2+4=6$ since the first block is not larger
than any other block, the second block is larger than the third and
the fifth blocks, and the third and the fourth blocks are both larger
than the fifth block.

The new statistic we now introduce is the sum of the statistic $\bmaj$
and the statistic $\mil$ computed on ordered partitions in the same
way as for unordered partitions.  (Thus, on ordered partitions, $\mil$
is \emm not, equal to $\los$ as it is on $\p n,k,$.)

\begin{definition}\label{bmm-def}%
  Let $\pi$ be an ordered partition with blocks $B_1,B_2,\ldots,B_k$
  and let $b_i$ be the size of block $i$.  Then
  
$$
\bmajmil \pi =  \bmaj\pi + \sum_i{(i-1)b_i}.
$$
\end{definition}  
For example, $\bmajmil(41-96-5-87-32) =
 (2+4) + (1\cd2+2\cd1+3\cd2+4\cd2)= 6+18=24$.

Note that $\mil$ is closely related to the $\maj$-statistic on
permutations.  Namely, $\maj$ of a permutation $p$ can be computed by
assigning to each letter in $p$ the number of descents to its right
and then summing these numbers.  This gives the same result as writing
$p$ backwards, cutting it at each non-descent and then computing
$\mil$ on the resulting partition.  As an example, 
\begin{eqnarray*}
\maj (2\,8\,6\,1\,3\,7\,4\,5) &=& 3+3+2+1+1+1+0+0\\[4pt]
\mil( 54-731-6-82) &=& 0+0+1+1+1+2+3+3.
\end{eqnarray*} 

In order to prove that $\bmajmil$ has the distribution $[k]!\sq n,k,$ we
need to rewrite identity \pref{zz}.  Recall that $\aq n,k,$ is the
polynomial in $q$ whose $i$-th coefficient is the number of
permutations in $\cls_n$ with $k$ descents and $\maj$ equal to $i$.

\begin{lemma}\label{aq-duality}%
  ~~$\disp \aq n,i-1, =q^{ni-{\ch n+1, 2,}}\cd \aq n,n-i,.  $
\end{lemma}
\begin{proof}{%
  If we reverse a permutation in $\cls_n$ with $i-1$ descents, the
  resulting permutation will have exactly $n-i$ descents.  Moreover,
  it is easy to see that the change in $\maj$ when reversing a
  permutation depends only on $n$ and $i$.  It is also easy to see
  that the polynomials $\aq n,k,$ are symmetric (reverse a permutation
  and complement each letter, i.e., replace $k$ by $n+1-k$).  Thus it
  suffices to compute the minimum of $\maj$ on the set of permutations
  with $i-1$ descents and on the set of those with $n-i$ descents,
  respectively, in order to determine the shift in the exponents of
  $q$ when going from $\aq n,i-1,$ to $\aq n,n-i,$.  The difference
  between these two minima is of course
\begin{eqnarray*}
\sum_{k=1}^{i-1}{i} - \sum_{k=1}^{n-i}{i} = \ch i,2,- \ch n-i+1,2,=
(i\cd(i-1) - (n-i+1)(n-i))/2,
\end{eqnarray*} 
which is easily seen to equal $ni-\ch n+1,2,$.
 }\end{proof}

Using Lemma~\ref{aq-duality}, we can now rewrite Proposition~\ref{zz}
in the following form.
\begin{proposition}[Equivalent to Proposition \protect{\ref{zz}}]\label{newzz}%
$$
[k]!\cd\sq n,k,=
\sum_i q^{k(k-i)+ni-{\ch n+1, 2,}}\cd\qbino n-i,k-i,\cd\aq n,n-i,.
$$
\end{proposition}

It is well-known that the $q$-binomial coefficient $\qbino n,k,$
records the inversion statistic on binary words of length $n$ with
exactly $k$ 1's.  As an example, $\qbino 3,2,=1+q+q^2$, corresponding
to the words $011,~ 101,~ 110$, whose number of inversions is zero,
one and two, respectively.  We will be concerned with inversions in
binary words associated to subsets of the set of descents of certain
permutations.

The \emm descent blocks, of a permutation $p$ are the maximal
contiguous decreasing subsequences of $p$.  For example, the
permutation $521364$ has descent blocks $521-3-64$.  

Let $p$ be a permutation in $\cls_n$ with descent blocks
$B_1,B_2,\ldots,B_i$.  Then $p$ has $n-i$ descents.  Pick $k-i$ of
these descents and mark them with 1's and mark the descents not picked
with 0's.  Let $w$ be the binary word obtained by reading the 0's and
1's from left to right.  We translate the permutation $p$, together
with the set of chosen descents (coded by $w$) into an ordered
partition $\pi$ with $k$ blocks by cutting $p$ into its descent blocks
and then further cutting the descent blocks at the chosen descents.

As an example, let $p= 51-742-6-83$, where we use the dashes to
indicate the partition of $p$ into descent blocks.  If we pick the
descents indicated by vertical bars in $5|1-74|2-6-83$ then $w=1010$
and the ordered partition obtained from $(p,w)$ is
$\pi=5-1-74-2-6-83$.  Now $\maj p=1+3+4+7=15$, $\inv w=3$ and
$\bmajmil( 5-1-74-2-6-83) = (1+4+3+4+10) + (1+3)=26$.  Thus,
$\bmajmil\pi=\maj p+\inv w +8$.  As it happens, the number 8 here
equals $k(k-i)+ni-{\ch n+1, 2,}$ (cf. Proposition~\ref{newzz}), since
we have $n=8$, $k=6$ (the number of blocks in $\pi$) and $i=4$ (the
number of descent blocks in $p$).

We will now prove that this holds in general, thus showing that the
distribution of $\bmajmil$ on $\op n,k,$ is given by the right hand
side in Proposition~\ref{newzz}.

\begin{theorem}%
  The distribution of $\bmajmil$ on $\op n,k,$ is 
$$
  \sum_i{q^{k(k-i)+ni-{\ch n+1, 2,}}\cd\qbino n-i,k-i,\cd\aq
    n,n-i,}.
$$
  Thus, $\bmajmil$ is Euler-Mahonian on ordered partitions.
\end{theorem}
\begin{proof}{%
  The proof is along the lines in the example preceding the theorem.
  That is, to each pair $(p,w)$, where $p$ is a permutation in
  $\cls_n$ with $n-i$ descents and $w$ is a binary word of length
  $n-i$ with $k-i$ 1's, we associate a partition $\pi$ by cutting $p$
  at each non-descent and at those descents in $p$ that correspond to
  1's in $w$.  This is clearly a bijective correspondence since each
  ordered partition in $\op n,k,$ gives rise to a unique permutation
  in $\cls_n$ with $n-i$ descents, exactly $k-i$ of which are marked
  (by virtue of occurring between adjacent blocks of $\pi$).
  
  We show that the difference $\bmajmil\pi -(\inv w + \maj p)$ equals
  the exponent $k(k-i)+ni-{\ch n+1, 2,}$ in the sum, which
  establishes the theorem.

  We fix $n$ and $i$, and proceed by induction on $k$.  The base case
  is $k=i$, since $k\ge i$.  The induction step consists of two parts,
  the first of which is a lemma that shows we may restrict our
  attention to words $w$ whose last letter is a 0.

\begin{itemize}
\item[(i)]%
  If $k=i$, referring to the example preceding the theorem, we have a
  permutation with no marked descents, so $w$ is the all zero word,
  whose $\inv$ is 0.  Thus we need to show that given a permutation
  $p$ with $k$ descent blocks, its $\maj$ differs from $\bmajmil$ of
  the corresponding ordered partition $\pi$ by the exponent
  $k(k-i)+ni-{\ch n+1, 2,}=0+n(k-1)-{\ch n, 2,} = nk-{\ch n+ 1, 2,}$.

Let $b_1, b_2,\cdots, b_k$ be the sizes of the descent blocks in $p$.
Then 
\begin{eqnarray*}
\maj p&=& {\ch n+ 1, 2,}- b_1 - (b_1 +b_2) -\cdots - (b_1 +b_2
+\cdots+ b_k)\\
&=& \ch n+ 1, 2, - \sum_{i=1}^k{(k+1-i)b_i}
\end{eqnarray*}
and, since $\bmaj\pi=0$, we have
$$\bmajmil\pi= b_2 +2b_3 +\cdots+ (k-1)b_k =
\sum_{i=1}^k{(i-1)b_i}.
$$  
Thus 
$$
\bmajmil\pi-\maj p = k\sum b_i -{\ch n+1, 2,} = kn - {\ch n+1, 2,},
$$
as desired.

\item[(ii)]%
  We now show that given a permutation $p$ and a binary word $w$
  indicating which of the descents in $p$ have been chosen, if we
  transpose a 0 in $w$ with an adjacent 1 to its right then the change
  in $\bmajmil$ of the corresponding partition is 1, which equals the
  change in $\inv w$.  (The transposition in $w$ does not, of course,
  affect $\maj p$.)
  
  We first consider the case where the 0 and the 1 belong to the same
  descent block.  Then $\bmajmil$ of the corresponding partition also
  increases by 1 because one letter gets moved from a block to the
  next one, whereas the number and position of block-descents is
  unchanged.
  
  If the 0 and 1 belong to different blocks then two letters get moved
  up one block, but the block descent associated to the 1 is reduced
  by one.  Thus the total increase in $\bmajmil$ is 1.

\item[(iii)]%
  Now we observe what happens when we increase $k$ by 1. This
  corresponds to changing one of the 0's in $w$ to a 1. By repeated
  use of (ii) we may assume that the last letter in $w$ is a 0 and
  that we are changing that 0 to a 1.
  
  The exponent $k(k-i)+ni-{\ch n+1, 2,}$ now increases by $2k+1-i$,
  but $\inv$ decreases by $k-i$, since one inversion is lost for each
  of the 1's we had.  Thus the total increase on the permutation side
  is $k+1$.
  
  On the partition side, suppose that the last $j$ blocks in $\pi$ are
  singletons (and not the block preceding them).  Then the 0 that is
  being changed to a 1 lies between the last two letters in block
  $k-j$ so the block number of the last letter in block $k-j$ is
  increased by 1 and the same is true of each of the $j$ letters in
  the trailing singleton blocks.  This contributes an increase of
  $j+1$ to $\bmajmil$.  Moreover, we introduce a new block descent in
  position $k-j$.  Thus the total increase in $\bmajmil$ is also
  $k+1$.
\end{itemize}
\par
\vspace{-2\blsk}
}\end{proof}

\section{Other (conjectured) Euler-Mahonian statistics}

We now introduce several new statistics on ordered partitions.  All of
these, like the statistic $\bmajmil$, are sums of one statistic defined
on the elements of the ordered partition and one statistic defined on
the partial ordering of the blocks described above.

All the conjectures in this section have been verified by computer for
all $n\le11$.  This is strong evidence, because what could ``go
wrong'' we expect to do so for much smaller $n$.

We start with a preliminary definition of a partition statistic that
we call $\mak$.  Namely, we let $\clos_i$ be the closer of the $i$-th
block in $\pi$ and set
$$
\mak\pi= \sum_i{\left(n-\clos_i\right)} +\rsb\pi. 
$$
This statistic is essentially the same as the Mahonian \emm
permutation, statistic bearing the same name, that was defined by
Foata and Zeilberger \cite{foata-zeilberger-denerts}.  Namely, let $\ncl\pi$ be the sum of
the non-closers in $\pi$ and observe that $\sum_i{\clos_i}\ = \ch
n+1,2,-\ncl\pi$.  Then $\mak\pi$ can be rewritten in the following
way, assuming $\pi$ has $k$ blocks:
\begin{eqnarray}\label{mak-re}%
\sum_i{\left(n-\clos_i\right)} +\rsb\pi
&=& kn-\sum_i{\clos_i} +\rsb\pi\\
&=& kn - \ch n+1,2,+\ncl\pi +\rsb\pi.\nonumber
\end{eqnarray} 
Now, for $n$ and $k$ fixed, the terms $kn - \ch n+1,2,$ only amount
to a constant shift of the statistic.  What remains is $\ncl\pi
+\rsb\pi$ which is equal to $\mak p$, where $p$ is the permutation
obtained from $\pi$ by writing the elements of the first block of
$\pi$ in decreasing order, then those of the second block in
decreasing order and so on, assuming that $\pi$ was written in
standard form.  As an example, the partition $421-63-85-7$ yields the
permutation $42163857$.  Thus, $\mak$ is, conjecturally, an
Euler-Mahonian statistic both on partitions and permutations.

There is another way to rewrite $\mak$ that is also interesting.
Namely, let $\scl_i$ be the number of closers in $\pi$ that are
smaller than the $i$-th letter in~$\pi$.  Then $\sum_i{(n-\clos_i)} =
\sum_i{\scl_i}$. Moreover, since no letter is larger than the closer
of its own block, we have that $\scl_i = \lcs_i + \rcs_i$.  Thus, as
we also have $\rsb = \ros-\rcs$, this entails that
\begin{eqnarray*}
\mak &=& \sum_i{\left(n-\clos_i\right)} +\rsb = 
\sum_i{\scl_i}  +\rsb \\[.5\blsk]
&=&  \lcs + \rcs + \rsb
= \lcs + \rcs + (\ros-\rcs) = \lcs  + \ros.\nonumber
\end{eqnarray*} 
On ordered partitions the statistic $\mak' = \lob+\rcb$ is clearly
equidistributed with $\mak$, which can be seen by complementing the
letters of a partition.  Moreover, by the same reasoning, the
statistics $n(k-1)-(\lcb+\rob)$ and $n(k-1)-(\los+\rcs)$ are
equidistributed with each other.
 
For clarity we now write up the definitions of these four statistics.

\begin{definition}
\label{mak-def}
\begin{eqnarray*}
\mak  &=& \lcs + \ros,\\
\mak' &=& \lob + \rcb,\\
\lmak &=& n(k-1) - (\los + \rcs),\\
\lmak' &=& n(k-1) - (\lcb + \rob).
\end{eqnarray*} 
\end{definition}

We conjecture that each of the above statistics can be paired with
either $\binv$ or $\bmaj$ to yield an Euler-Mahonian statistic on
ordered partitions.

\begin{conjecture}\label{omak-conj}%
  Each of the following eight statistics is Euler-Mahonian on $\op
  n,k,$.
\begin{eqnarray*}
&&\mak + \bmaj, ~~~~  \mak' + \bmaj, ~~~~ \lmak  + \bmaj, ~~~~ \lmak'
+ \bmaj,\\ 
&&\mak + \binv, ~~~~~ \mak' + \binv, ~~~~\, \lmak  + \binv, ~~~~\, \lmak'  +
\binv.
\end{eqnarray*}
\end{conjecture}
Since both $\binv$ and $\bmaj$ vanish on $\p n,k,$ the above eight
statistics reduce to the statistics in Definition~\ref{mak-def} when
restricted to $\p n,k,.$ It was pointed out by an anonymous referee
that they are all Euler-Mahonian on $\p n,k,$, that $\mak'$ and
$\lmak$ are both equal (as functions) to $\RCB$ (see
Theorem~\ref{ros-thm}), and that $\mak$ and $\lmak'$ are equal, and
can be shown to be Euler-Mahonian using a bijection in Section~3 of
\cite{white-interpolating}.

The statistics in Definition~\ref{mak-def} were also shown to be
Euler-Mahonian by G\'erald Ksavrelof and Jiang Zeng
\cite{ksavrelof-zeng-nouvelle-stat}.  In fact, they proved more.  Namely, given any
partition $\pi=B_1-B_2-\cdots -B_k$ of $[n]$ into $k$ blocks and $d$
with $0\leq d\leq k$, define
$$
\mak_d\pi=\mak\pi+k-d-\#\{a\in B_j\st j>d,\; a>\clos_d\}.
$$
Then they prove that $\mak_d$ is Euler-Mahonian on $P_n^k$. Note
that $\mak_k=\mak$.

\smallskip

The last new statistic we introduce is not based on any Mahonian
permutation statistic, but rather on the coordinate statistic $\lsb$
introduced in Definition~\ref{stats-def}.  We could of course make an
analogous definition for $\rsb$, which would require a ``right''
analog of $\cbmaj$ and $\cbinv$.

\begin{conjecture}\label{lsb-conj}%
The statistics 
\begin{eqnarray*}
\cmajlsb &=& \lsb + \cbmaj + \ch k,2,,\\
\cinvlsb &=& \lsb + \cbinv + \ch k,2,
\end{eqnarray*} 
are Euler-Mahonian on \op n,k,.
\end{conjecture}

Note that when we restrict to $\p n,k,$, both $\cmajlsb$ and
$\cinvlsb$ reduce to $\lsb +\ch k,2,$.  On $\p n,k,$ the statistic
$\lsb$ is equal to $\lcb$, because each block to the left of the block
containing a given letter $i$ has an opener that is smaller than $i$.
Thus, by Corollary~\ref{etc-coro},  $\lsb +\ch k,2,$  is
Euler-Mahonian on $\p n,k,$.

Finally we point out an equivalent way of formulating the conjecture
that $\mak+\bmaj$ and $\mak+\binv$ are Euler-Mahonian.  This
formulation gives a simpler form for the sum in
Proposition~\ref{newzz} because the shift in $\mak$ mentioned above
can be factored out from that sum in a way that allows us to eliminate
the power of $q$ involved.  We use the following lemma, where $\qbino
a,b,_{q^{-1}}$ is the polynomial in $q^{-1}$ obtained by replacing $q$
with $q^{-1}$ in $\qbino a,b,$:
\begin{lemma}\label{qbino-lemma} %
$$
\qbino n-i, k-i,_{q^{-1}} = q^{(k-n)(k-i)}\qbino n-i, k-i,.
$$
\end{lemma}
\begin{proof}{%
  It is well known that $\qbino a,b,$ is symmetric as a polynomial in
  $q$.  Therefore, since the effect of replacing $q$ by $q^{-1}$ in
  the polynomial $\qbino a,b,$ is to reverse the coefficients and
  shift the exponents of $q$ in $\qbino a,b,$, to prove the lemma we
  only need to compute the degree of $\qbino n-i, k-i,$ as a
  polynomial in $q$.  This degree is $(n-k)(k-i)$.
}\end{proof}

Let now $M_q(n,k)$ be the distribution of $\ncl + \rsb +\bmaj$ on $\op
n,k,$.  Then, by identity~\pref{mak-re}, $M_q(n,k)$ equals the
distribution  on $\op n,k,$ of $\mak+\bmaj$, divided by $q^{kn-\ch
  n+1,2,}$. Thus, the conjecture that $\mak+\bmaj$ is Euler-Mahonian
on $\op n,k,$ is equivalent to 
\begin{eqnarray*}%
  M_q(n,k) = q^{\ch n+1,2,-kn}\sum_i q^{k(k-i)+ni-{\ch n+1,
      2,}}\cd\qbino n-i,k-i,\cd\aq n,n-i,.
\end{eqnarray*} 
Rewriting the power of $q$ in the right hand side and then applying
Lemma~\ref{qbino-lemma} we obtain the following, which is equivalent
to the first part of Conjecture~\ref{omak-conj}.
\begin{conjecture} \label{strip-conj}%
\begin{eqnarray*}%
  M_q(n,k) = \sum_i{\qbino n-i,k-i,_{q^{-1}}\cd\aq n,n-i,}.
\end{eqnarray*} 
\end{conjecture}
Observe that we could also take $M_q(n,k)$ to be the distribution of
$\ncl + \rsb +\binv$.  What is perhaps most interesting about
Conjecture~\ref{strip-conj} is the fact that $\aq n,n-i,$ is also the
distribution of the permutation statistic $\mak$ on the set of
permutations in $\cls_n$ with exactly $n-i$ descents.  (In other
words, the bistatistics $(\des,\mak)$ and $(\des,\maj)$ have the same
distribution on $\cls_n$.)  One might thus hope to find a simple proof
of Conjecture~\ref{strip-conj}, since we are computing essentially the
same statistic on both sides, but we have been unable to find such a
proof.

\bigskip

\centerline{\sc Acknowledgments}

\smallskip

My thanks to Jiang Zeng for fruitful discussions in the initial stage
of the research presented here, and to the anonymous referee mentioned
after Conjecture~\ref{omak-conj}.

\bibliographystyle{abbrv}
\bibliography{allrefs}

\end{document}